\documentclass[12pt]{article}
\usepackage[utf8]{inputenc}
\usepackage{amssymb,amsfonts}
\usepackage[english,russian]{babel}
\usepackage[usenames]{color}
\usepackage{wrapfig}
\usepackage{url}
\usepackage{amsmath}
\usepackage{cmap}
\usepackage{graphicx}
\usepackage{graphics}

\newcommand{\2}{<<}
\newcommand{\3}{>>}

\newcommand{\8}{$\blacktriangleleft\,$}
\newcommand{\9}{$\,\blacktriangleright$}
\newcommand{\q}{\mathbb{Q}}

\newcommand{\R}{\mathbb{R}}

\newcommand{\fo}[1]{${\displaystyle #1}$}

\newcommand{\mr}{\leqslant}
\newcommand{\br}{\geqslant}
\newcounter{teorema}
\newcounter{glav}
\newcounter{risunok}
\newcommand{\ris}{\refstepcounter{risunok}{\bf Рис. \arabic{risunok}.}}
\newcommand{\glava}[1]{\refstepcounter{glav}\begin{center}\textbf{\arabic{glav}. #1}\end{center}}

\newcommand{\teop}{\refstepcounter{teorema}{\bf Теорема \arabic{teorema} }}
\newcommand{\problem}{\refstepcounter{prob}\textbf{Задача \arabic{prob}. }}

\renewcommand{\textrm}[1]{\textit{#1}}
\newcounter{le}
\newcounter{op}
\newcounter{prob}
\newcounter{exa}
\newcommand{\lemma}{\refstepcounter{le}\textbf{Лемма \arabic{le}. }}
\newcommand{\lemmap}{\refstepcounter{le}\textbf{Лемма \arabic{le} }}
\newcommand{\opred}{\refstepcounter{op}\textbf{Определение \arabic{op}. }}
\begin{document}
\begin{center}
{\Large \textbf{\fo{\textit{x}}-площадь}}\\
$\phantom{kl}$\\
\begin{large}
\textbf{Фёдор Шаров\footnote{Автор был частично поддержан грантом президента Российской Федерации MK-6137.2016.1.}}\\
$\phantom{kl}$\\
\end{large}
{\small\textit{Факультет математики, Национальный исследовательский университет Высшая школа экономики}}
\end{center}

С давних времён люди занимались задачами на разрезания. Можно ли разрезать определённую фигуру на какие-то другие (например, на квадраты)? Можно ли из некоторого набора фигур составить другую заданную фигуру (например, правильный треугольник)? Такими вопросами люди задавались с незапамятных времён. В этой статье мы докажем один из классических результатов в этой области математики --- теорему Дена.

\teop\label{D}\textbf{ (М.~Ден, 1903).} \textit{Если прямоугольник разрезан на квадраты (не обязательно равные), то отношение его сторон рационально.}

%У этой теоремы существует множество доказательств. Первое из них было получено М.~Деном в 1903~г. и было довольно сложным. Впоследствии, доказательство несколько раз упрощалось. Например, есть красивый метод доказательства, принадлежащий Р.~Л.~Бруксу, К.~А.~Б.~Смиту, А.~Г.~Стоуну и У.~Т.~Татту, в котором разрезанию ставится в соответствие электрическая цепь~\cite{KV1}. Другой метод доказательства~\cite{yagl}, принадлежащий И. М. Яглому, заключается в очень маленьком увеличении или уменьшении сторон разрезаемого прямоугольника и квадратов разрезания, а затем доказывается, что непрерывно такую операцию производить нельзя --- разрезание \2разрушится\3. Наконец, есть неэлементарный метод, использующий аддитивные функции.

%В §II.1 приведён набор задач, решение которых позволяет доказать теорему Дена при помощи метода \2площадей\3. В §II.2 приведены решения наиболее интересных задач.

%\begin{center}\textbf{§II.1. Задачи}\end{center}

\begin{wrapfigure}[9]{l}{0.26\linewidth} 
\vspace{-3ex}
\begin{picture}(45,45) 
\put(5,43){\line(1,0){15}} 
\put(5,37){\line(1,0){15}}
\put(5,31){\line(1,0){15}}
\put(26,43){\line(1,0){15}} 
\put(26,37){\line(1,0){15}}
\put(26,31){\line(1,0){15}}
\put(5,5){\line(1,0){15}} 
\put(5,11){\line(1,0){15}}
\put(5,17){\line(1,0){15}}
\put(26,5){\line(1,0){15}} 
\put(26,11){\line(1,0){15}}
\put(26,17){\line(1,0){15}}
\put(5,43){\line(0,-1){15}}
\put(11,43){\line(0,-1){15}}
\put(17,43){\line(0,-1){15}}
\put(29,43){\line(0,-1){15}}
\put(35,43){\line(0,-1){15}}
\put(41,43){\line(0,-1){15}}
\put(5,20){\line(0,-1){15}}
\put(11,20){\line(0,-1){15}}
\put(17,20){\line(0,-1){15}}
\put(29,20){\line(0,-1){15}}
\put(35,20){\line(0,-1){15}}
\put(41,20){\line(0,-1){15}}
\put(21,43){\circle*{0.5}}
\put(23,43){\circle*{0.5}}
\put(25,43){\circle*{0.5}}
\put(21,5){\circle*{0.5}}
\put(23,5){\circle*{0.5}}
\put(25,5){\circle*{0.5}}
\put(5,22){\circle*{0.5}}
\put(5,24){\circle*{0.5}}
\put(5,26){\circle*{0.5}}
\put(41,22){\circle*{0.5}}
\put(41,24){\circle*{0.5}}
\put(41,26){\circle*{0.5}}
\put(3.5,5.5){$^x$}
\put(3.5,11.5){$^x$}
\put(3.5,31.5){$^x$}
\put(3.5,37.5){$^x$}
\put(7.5,42){$^x$}
\put(13.5,42){$^x$}
\put(37.5,42){$^x$}
\put(31.5,42){$^x$}
\put(20,2.1){$mx$}
\put(41.5,23.5){$nx$}
\put(17,-4){\ris\label{ri2}}
\end{picture}
\end{wrapfigure}

Читатель, конечно, без труда сможет догадаться, как разрезать прямоугольник с рациональным отношением сторон на квадраты (см. рис.~\ref{ri2}). Содержательным является лишь вопрос невозможности разрезаний в остальных случаях.

%Сперва, предлагаем читателю самому попробовать доказать теорему Дена, решив задачи, приведённые ниже. Мы попробуем сначала подвести читателя к доказательству теоремы для частного случая, когда стороны разрезаемого прямоугольника и стороны всех квадратов являются числами вида \fo{a+b\sqrt{2}}, где $a$ и $b$ рациональны. А затем, если читатель справится с этими задачами, ему будет предложено попробовать доказать теорему Дена и в общем случае. Дальнейший текст статьи содержит решения наиболее важных из предложенных задач, а также описание других результатов, полученных применением того же метода, что используется в нашем доказательстве теоремы Дена.

\textbf{Соглашения.} Латинские буквы $a$, $b$, $c$, $d$ и эти же буквы с индексами в этой статье означают \textit{рациональные} числа. Все числа, которые можно представить в виде ${x=a+b\sqrt{2}}$, будем называть \textit{хорошими}.

\glava{Теорема Дена для хороших чисел}

Сначала докажем теорему для частного случая, когда стороны разрезаемого прямоугольника и стороны всех квадратов являются хорошими числами.

Прежде всего рассмотрим пример: прямоугольник \fo{1\times\sqrt{2}}. Решим следующую задачу.

\problem\label{z1} Можно ли прямоугольник \fo{1\times\sqrt{2}} разрезать на квадраты с рациональными сторонами? А со сторонами, которые либо рациональны, либо имеют вид \fo{b\sqrt{2}}? А со сторонами, которые являются произвольными хорошими числами?

\textbf{Решение (написано Л.~Алиевой).} \textit{Ответ:} на первые два вопроса --- нет. Решение для третьего вопроса см. в конце статьи.

\8Ответим на первый вопрос. Предположим, что прямоугольник \fo{1\times\sqrt{2}} разрезан на квадраты с рациональными сторонами. Рассмотрим сторону, равную \fo{\sqrt{2}}. Сумма примыкающих к ней сторон квадратов должна равняться $\sqrt{2}$. Сумма рациональных чисел не может быть равна иррациональному, следовательно, такое разрезание невозможно.

Ответим на второй вопрос. Предположим, что прямоугольник \fo{1\times\sqrt{2}} разрезан на квадраты со сторонами вида $a$ или $b\sqrt{2}$. Его площадь равна $\sqrt{2}$. Площади квадратов имеют вид либо $a^2$, либо \fo{2b^2}. Площадь прямоугольника равна сумме площадей квадратов. Получаем противоречие, так как иррациональное число не может быть равно сумме рациональных. Здесь мы опирались на \textit{аддитивность} площади: площадь целого равна сумме площадей частей.\9

%\problem\label{z1.5} Ответьте на вопросы задачи~\ref{z1} для прямоугольников \fo{1\times\left(1+\sqrt{2}\right)} и \fo{1\times\left(2+\sqrt{2}\right)}.
%Попробуем теперь ответить на вопросы задачи~\ref{z1} для прямоугольников \fo{1\times\left(1+\sqrt{2}\right)} и \fo{1\times\left(2+\sqrt{2}\right)}.

Покажем теперь, что прямоугольник \fo{1\times\left(1+\sqrt{2}\right)} нельзя разрезать на квадраты с хорошими сторонами. Для этого нам понадобится определение \textit{сопряжённого числа}.

\opred Число \fo{\overline{s}:=a-b\sqrt{2}} назовём \textit{сопряжённым} к числу \fo{s=a+b\sqrt{2}}.

\problem\label{z_sopr} Докажите, что сопряжённое к сумме хороших чисел равно сумме сопряжённых, а сопряжённое к произведению --- произведению сопряжённых.

\footnote{Следующий текст с доказательством невозможности разрезания прямоугольника \fo{1\times\left(1+\sqrt{2}\right)} на квадраты написан Л.~Алиевой.}\8Итак, предположим, что прямоугольник \fo{1\times\left(1+\sqrt{2}\right)} разрезан на $n$ квадратов со сторонами \fo{a_1+b_1\sqrt{2}}, \fo{a_2+b_2\sqrt{2}}, \ldots, \fo{a_n+b_n\sqrt{2}}. Площадь данного прямоугольника равна сумме площадей рассматриваемых $n$~квадратов, то есть $$1+\sqrt{2}=\left(a_1+b_1\sqrt{2}\right)^2+\left(a_2+b_2\sqrt{2}\right)^2+\ldots+\left(a_n+b_n\sqrt{2}\right)^2.$$ Найдём сопряжённые к обоим частям полученного уравнения. Сопряжённым к числу \fo{1+\sqrt{2}} является число \fo{1-\sqrt{2}}. По утверждению задачи~\ref{z_sopr} получаем, что сопряжённым к числу $$\left(a_1+b_1\sqrt{2}\right)^2+\left(a_2+b_2\sqrt{2}\right)^2+\ldots+\left(a_n+b_n\sqrt{2}\right)^2$$ является число $$\left(a_1-b_1\sqrt{2}\right)^2+\left(a_2-b_2\sqrt{2}\right)^2+\ldots+\left(a_n-b_n\sqrt{2}\right)^2.$$ Получаем, что $$1-\sqrt{2}=\left(a_1-b_1\sqrt{2}\right)^2+\left(a_2-b_2\sqrt{2}\right)^2+\ldots+\left(a_n-b_n\sqrt{2}\right)^2.$$ Заметим, что \fo{1-\sqrt{2}} --- число отрицательное, а правая часть тождества --- сумма неотрицательных. Получаем противоречие. Значит, прямоугольник \fo{1\times\left(1+\sqrt{2}\right)} нельзя разрезать на квадраты с хорошими сторонами.\9

%Теперь покажем, что прямоугольник \fo{1\times\left(2+\sqrt{2}\right)} нельзя разрезать на квадраты с хорошими сторонами.

Заметим, что, например, сопряжённое к числу \fo{2+\sqrt{2}} положительно. Поэтому только что использованный метод доказательства невозможности разрезания прямоугольника \fo{1\times\left(1+\sqrt{2}\right)} неприменим для доказательства невозможности разрезания прямоугольника \fo{1\times\left(2+\sqrt{2}\right)}, а его тоже нельзя разрезать на квадраты с хорошими сторонами согласно теореме Дена. Обойти это препятствие можно при помощи следующего обобщения понятия площади.% И доказать, что эта новая \2площадь\3 также удовлетворяет свойству аддитивности (лемма~\ref{pr4} ниже).% Но выход есть. Подумаем, какими свойствами площади мы пользовались при доказательстве невозможности разрезания прямоугольника \fo{1\times\left(1+\sqrt{2}\right)}? Мы вновь пользовались свойством аддитивности площади. А нельзя ли придумать какую-нибудь функцию, которая каждому прямоугольнику будет сопоставлять какое-то число так, чтобы для этой функции выполнялось свойство аддитивности? Оказывается можно.

\opred Пусть $x$ --- некоторое действительное число. Назовём \textit{\fo{x\mbox{-площадью}}} (или \textit{площадью Гамеля}) прямоугольника ${\left(a+b\sqrt{2}\right)\times\left(c+d\sqrt{2}\right)}$ число \fo{(a+bx)(c+dx)}.

%Вряд ли получится ответить на вопросы задачи~\ref{z1.5} для прямоугольника \fo{1\times\left(2+\sqrt{2}\right)} без следующего обобщения понятия площади (мы обобщаем понятие площади так, чтобы площадь этого прямоугольника стала отрицательной, а площади квадратов оставались неотрицательными).

%Каждый из пунктов а)--г) можно решать, основываясь на невозможности равенства площади разрезаемого прямоугольника сумме площадей квадратов, на которые производится разрезание. При этом используется свойство площади под названием \textit{аддитивность}: сумма площадей фигур, на которые мы разрезаем, должна быть равна площади разрезаемой фигуры. Попробуем обобщить понятие площади так, чтобы площадь прямоугольника \fo{1\times\left(2+\sqrt{2}\right)} стала отрицательной, а площади квадратов оставались неотрицательными. Пункт д) вряд ли получится решить без следующего понятия и его свойств (задачи~\ref{pr3} и~\ref{pr4}). \textbf{\textit{(Не кажется ли Вам, что в этом абзаце первые два предложения слишком сильно помогают решить пункты а)--г)? Или это нормально для первой задачи?)}}

\problem Докажите, что обычная площадь прямоугольника ${\left(a+b\sqrt{2}\right)\times\left(c+d\sqrt{2}\right)}$ и сопряжённое к ней число --- это одни из его \fo{x\mbox{-площадей}}. Чему равно $x$ в каждом из случаев?

%\problem Докажите, что сопряжённое число к обычной площади прямоугольника --- одна из его \fo{x\mbox{-площадей}}. Чему равно $x$ в этом случае?

\problem\label{z4} Найдите все прямоугольники вида ${\left(a+b\sqrt{2}\right)\times\left(c+d\sqrt{2}\right)}$, \fo{x\mbox{-площади}} которых неотрицательны при всех~$x$.

%\problem Докажите, что для любого прямоугольника ${\left(a+b\sqrt{2}\right)\times\left(c+d\sqrt{2}\right)}$, отношение сторон которого иррационально, существует такое число~$x$, что \fo{x\mbox{-площадь}} этого прямоугольника отрицательна.

%Заметим ещё, что $x$-площадь может быть и отрицательным числом. Но в тех случаях, когда числа \fo{a+bx} и \fo{c+dx} положительны, $x$-площадь прямоугольника ${\left(a+b\sqrt{2}\right)\times\left(c+d\sqrt{2}\right)}$ имеет простой геометрический смысл: она равна обычной площади прямоугольника ${\left(a+bx\right)\times\left(c+dx\right)}$.

%\refstepcounter{prob}\textbf{Задача \arabic{prob}* (не используется в дальнейшем).} Пусть прямоугольник ${1\times\left(A+B\sqrt{2}\right)}$ можно разрезать на квадраты с длинами сторон вида \fo{a+b\sqrt{2}}. Заменим в выражении для длины каждого квадрата число~\fo{\sqrt{2}} на  число~$x$, которое выберем таким образом, чтобы длина стороны любого из этих квадратов оставалась положительной. Верно ли что из полученных квадратов можно сложить прямоугольник ${1\times\left(A+Bx\right)}$?

%\refstepcounter{prob}\textbf{Задача \arabic{prob} (неотрицательность \textit{x}-площади квадрата).}\label{pr3} Докажите, что для любого $x$ $x$-площадь квадрата ${\left(a+b\sqrt{2}\right)\times\left(a+b\sqrt{2}\right)}$ неотрицательна.

\lemmap\textbf{(аддитивность \textit{x}-площади).}\label{pr4} \textit{Если прямоугольник разрезан на конечное число прямоугольников, стороны которых --- хорошие числа, то для любого ${x\in\R}$ ${x\mbox{-площадь}}$ разрезаемого прямоугольника равна сумме ${x\mbox{-площадей}}$ прямоугольников, на которые он разрезан.}
%произвольное число прямоугольников со сторонами вида \fo{a+b\sqrt{2}}, то для любого $x$ \fo{x\mbox{-площадь}} разрезаемого прямоугольника равна сумме \fo{x\mbox{-площадей}} прямоугольников, на которые он разрезан.

%Предлагаемый метод доказательства невозможности разрезаний состоит в том, чтобы найти такое $x$, что \fo{x\mbox{-площадь}} разрезаемой фигуры была бы отрицательной, а \fo{x\mbox{-площадь}} любой из фигур, на которые мы разрезаем --- неотрицательной. Далее из аддитивности \fo{x\mbox{-площади}} (лемма~\ref{pr4}) заключаем, что разрезание невозможно.

\begin{wrapfigure}[10]{l}{0.28\linewidth}
\vspace{-2ex}%-2
\hspace*{6ex}
\begin{picture}(45,83)
\thicklines
\put(0,80){\line(1,0){35}}
\put(0,60){\line(1,0){35}}
\put(0,60){\line(0,1){20}}
\put(35,60){\line(0,1){20}}
\put(14,60){\line(0,1){20}}
\put(5.5,68){$S_1$}
\put(22.5,68){$S_2$}
\put(-11,68){$^{a+b\sqrt{2}}$}
\put(0.5,80){$^{c_1+d_1\sqrt{2}}$}
\put(18.5,80){$^{c_2+d_2\sqrt{2}}$}
%\put(48,80){$^{\gamma+\delta\sqrt{p}}$}

%\put(42,80){\line(1,0){21}}
%\put(42,60){\line(1,0){21}}
%\put(42,60){\line(0,1){20}}
%\put(63,60){\line(0,1){20}}
%\put(42,72){\line(1,0){21}}
%\put(50.5,74.5){$S_3$}
%\put(50.5,64.5){$S_4$}
%\put(63.5,74){$^{\alpha_1+\beta_1\sqrt{p}}$}
%\put(63.5,64){$^{\alpha_2+\beta_2\sqrt{p}}$}

\put(7.5,54){\ris\label{rrr6}}
\end{picture}
\end{wrapfigure}
 
%\textit{принадлежащими множеству \fo{\hor{p}}, то \2площадь\3 разрезаемого прямоугольника равна сумме \2площадей\3 прямоугольников, на которые он разрезан.}
\textbf{Доказательство леммы~\ref{pr4}.} $\blacktriangleleft\,$Нетрудно убедиться, что сумма \fo{x\mbox{-площадей}} двух прямоугольников со сторонами вида \fo{a+b\sqrt{2}} равна \fo{x\mbox{-площади}} их объединения. Действительно, пусть имеется прямоугольник с \fo{x\mbox{-площадью}}~\fo{S}, который состоит из двух прямоугольников с \fo{x\mbox{-площадями}}~\fo{S_1} и~\fo{S_2} (см. рис.~\ref{rrr6}). Тогда $$S_1+S_2=(a+bx)(c_1+d_1x)+(a+bx)(c_2+d_2x)=$$ $$
=(a+bx)((c_1+c_2)+(d_1+d_2)x)=S.$$

\begin{wrapfigure}[6]{l}{0.28\linewidth}%по умолчанию 6 строк
\vspace{-1.2ex}
\hspace*{2.5ex}
\begin{picture}(45,83)
\thicklines
\put(-5,75.3){$^{III}$}
\put(-3.8,71.3){$^{II}$}
\put(-2.6,63.3){$^I$}
\put(0,80){\line(1,0){45}}
\put(14,72){\line(1,0){31}}
\put(21,76){\line(1,0){24}}
\put(0,60){\line(1,0){45}}
\put(0,60){\line(0,1){20}}
\put(14,60){\line(0,1){20}}
\put(21,72){\line(0,1){8}}
\put(25,60){\line(0,1){12}}
\put(35,60){\line(0,1){12}}
\put(45,60){\line(0,1){20}}
\thinlines
\multiput(0,72)(2,0){8}{\line(1,0){1}}
\multiput(0,76)(2,0){11}{\line(1,0){1}}
\multiput(21,60)(0,2){6}{\line(0,1){1}}
\multiput(25,80)(0,-2){5}{\line(0,-1){1}}
\multiput(35,80)(0,-2){5}{\line(0,-1){1}}
\put(15.5,55){\ris\label{rrr5}}
\end{picture}
\end{wrapfigure}

Пусть теперь количество прямоугольников в разрезании больше двух. Продолжим каждый разрез, как показано на рис.~\ref{rrr5}. Тогда каждый прямоугольник нового разрезания будет также иметь стороны вида \fo{a+b\sqrt{2}}. Рассмотрим горизонтальные слои из последовательно приложенных друг к другу по общей стороне прямоугольников ($I$, $II$, $III$ на рис.~\ref{rrr5}). Используя уже доказанное свойство аддитивности \fo{x\mbox{-площади}} двух прямоугольников с общей стороной, методом математической индукции легко доказывается, что \fo{x\mbox{-площадь}} любого такого слоя равна сумме \fo{x\mbox{-площадей}} прямоугольников, составляющих этот слой. Теперь уже эти слои приложим друг к другу и применим только что доказанное утверждение об аддитивности \fo{x\mbox{-площади}} ряда прямоугольников. Получим, что \fo{x\mbox{-площадь}} разрезаемого прямоугольника равна сумме \fo{x\mbox{-площадей}} горизонтальных слоёв. Эта сумма равна сумме \fo{x\mbox{-площадей}} всех прямоугольников разрезания.$\,\blacktriangleright$

%\textit{Указание.} Начните со случая разрезания на два прямоугольника.

\teop\textbf{(частный случай теоремы Дена).}\label{std} \textit{Прямоугольник с хорошими длинами сторон можно разрезать на квадраты с хорошими длинами сторон тогда и только тогда, когда отношение его сторон рационально.}

\textbf{Доказательство.} $\blacktriangleleft\,$Достаточно доказать невозможность разрезания на квадраты прямоугольников вида \fo{1\times\left(c+d\sqrt{2}\right)}, где \fo{d\neq0}. Пусть такой прямоугольник разрезан на квадраты со сторонами вида \fo{a+b\sqrt{2}}. Его \fo{x\mbox{-площадь}} равна \fo{c+dx}. Так как \fo{d\neq0}, то эта \fo{x\mbox{-площадь}} отрицательна при некотором~$x$. В то же время, \fo{x\mbox{-площадь}} любого квадрата со стороной \fo{a+b\sqrt{2}} равна \fo{(a+bx)^2}, что неотрицательно для любого~$x$. Получили противоречие с леммой~\ref{pr4}: сумма неотрицательных чисел не может равняться отрицательному.$\,\blacktriangleright$

\glava{Доказательство теоремы Дена в общем случае}

Для доказательства теоремы Дена в общем случае определение \fo{x\mbox{-площади}} нам уже не годится: ведь она определена только для хороших чисел, а теперь у нас в разрезании могут присутствовать квадраты с какими угодно сторонами (например, \fo{\sqrt{3}}, или \fo{\pi}, или \fo{\sqrt[3]{2}} и так далее).

Далее в этом разделе мы считаем, что прямоугольник ${s_0\times t_0}$ разрезан на прямоугольники ${s_1\times t_1}$, ${s_2\times t_2}$, \ldots, ${s_N\times t_N}$, причём $s_0$ и $t_0$ несоизмеримы.

\lemma\label{l0} \textit{Обозначим
$$
P=\{s_0,t_0,s_1,t_1,\ldots,s_N,t_N\}.
$$
Тогда можно выбрать такие числа ${e_1,e_2,\ldots,e_n\in
P}$, чтобы любое число ${p\in P}$ единственным образом представлялось
в виде
$$
p=as_0+bt_0+a_1e_1+a_2e_2+\ldots+a_ne_n.
$$}

\begin{wrapfigure}[6]{l}{0.31\linewidth}
\vspace{-4.5ex}
\hspace*{1ex}
\begin{picture}(45,83)
\thicklines
%\put(0,80){\line(1,0){35}}
%\put(0,60){\line(1,0){35}}
\put(0,60){\line(0,1){15}}
\put(25,60){\line(0,1){15}}
\put(51,60){\line(0,1){15}}
\put(0,75){\line(1,0){51}}
\put(0,60){\line(1,0){51}}
\put(0,70){\line(1,0){25}}
%\put(35,60){\line(0,1){20}}
%\put(14,60){\line(0,1){20}}
%\put(5.5,68){$S_1$}
%\put(22.5,68){$S_2$}
\put(-2,66){1}
\put(20,75.5){\fo{2+\sqrt{2}}}
\put(7.5,69.8){\fo{^{1/3\times\sqrt{3}}}}
\put(7.5,62.3){\fo{^{2/3\times\sqrt{3}}}}
\put(27.8,64.8){\fo{^{1\times\left(2+\sqrt{2}-\sqrt{3}\right)}}}
%\put(0.5,80){$^{c_1+d_1\sqrt{2}}$}
%\put(18.5,80){$^{c_2+d_2\sqrt{2}}$}
%\put(48,80){$^{\gamma+\delta\sqrt{p}}$}

%\put(42,80){\line(1,0){21}}
%\put(42,60){\line(1,0){21}}
%\put(42,60){\line(0,1){20}}
%\put(63,60){\line(0,1){20}}
%\put(42,72){\line(1,0){21}}
%\put(50.5,74.5){$S_3$}
%\put(50.5,64.5){$S_4$}
%\put(63.5,74){$^{\alpha_1+\beta_1\sqrt{p}}$}
%\put(63.5,64){$^{\alpha_2+\beta_2\sqrt{p}}$}

\put(17.5,55){\ris\label{ri3}}
\end{picture}
\end{wrapfigure}

Например, для разрезания на рис.~\ref{ri3}, где \fo{s_0=1} и \fo{t_0=2+\sqrt{2}}, необходимо выбрать лишь одно число (либо \fo{e_1=\sqrt{3}}, либо \fo{e_1=2+\sqrt{2}-\sqrt{3}}).

\textbf{Доказательство.} \8Выпишем в строку длины сторон прямоугольников разрезания, начиная с $s_0$ и $t_0$. Например, для разрезания на рис.~\ref{ri3} получится такая строка:
%Попробуем лучше разобраться, что такое базис, на простом примере. Пусть у нас есть 9 чисел:
%$$1,\,\,\,\,\,z=2+\frac{\sqrt{2}}{4}+\frac{\sqrt{3}}{4},\,\,\,\,\,s_1=\frac{1}{2},\,\,\,\,\,s_2=\frac{\sqrt{2}}{2},\,\,\,\,\,s_3=\frac{1}{2}+\frac{\sqrt{3}}{4},\,\,\,\,\,s_4=\frac{1}{4}+\frac{\sqrt{2}}{2},\,\,\,\,\,s_5=\frac{1}{4},\,\,\,\,\,s_6=\frac{\sqrt{3}}{4}.$$
$$s_0=1,\,\,\,\,\,t_0=2+\sqrt{2},\,\,\,\,\,s_1=1/3,\,\,\,\,\,t_1=\sqrt{3},\,\,\,\,\,s_2=2/3,\,\,\,\,\,t_2=\sqrt{3},\,\,\,\,\,s_3=1,\,\,\,\,\,t_3=2+\sqrt{2}-\sqrt{3}.$$ \footnote{Дальнейший текст доказательства леммы~\ref{l0} с незначительными изменениями заимствован из статьи Д. Фукса \2Можно ли из тетраэдра сделать куб?\3 (Квант~№11 (1990), стр. 2--11).}Теперь подчеркнём в этом списке те числа, которые не представляются в виде суммы предыдущих с рациональными коэффициентами (сейчас поясним, что это значит). Например, $s_0$ мы подчеркнём (предыдущих чисел вовсе нет). Число $t_0$ мы тоже подчеркнём, так как оно несоизмеримо с $s_0$ и, следовательно, не представимо в виде ${a_0s_0}$. Далее, если $s_1$ представляется в виде ${a_0s_0+b_0t_0}$, то $s_1$ не подчёркиваем, а если не представляется, то подчёркиваем. Аналогично, если $t_1$ представляется в виде ${a_0s_0+b_0t_0+a_1s_1}$, то $t_1$ не подчёркиваем, а если не представляется, то подчёркиваем. И так далее. Например, в нашем случае получится такое: 
%$$\underline{1},\,\,\,\,\,\underline{2+\frac{\sqrt{2}}{4}+\frac{\sqrt{3}}{4}},\,\,\,\,\,\frac{1}{2},\,\,\,\,\,\underline{\frac{\sqrt{2}}{2}},\,\,\,\,\,\frac{1}{2}+\frac{\sqrt{3}}{4},\,\,\,\,\,\frac{1}{4}+\frac{\sqrt{2}}{2},\,\,\,\,\,\frac{1}{4},\,\,\,\,\,\frac{\sqrt{3}}{4}.$$
$$\underline{1},\,\,\,\,\,\underline{2+\sqrt{2}},\,\,\,\,\,1/3,\,\,\,\,\,\underline{\sqrt{3}},\,\,\,\,\,2/3,\,\,\,\,\,\sqrt{3},\,\,\,\,\,1,\,\,\,\,\,2+\sqrt{2}-\sqrt{3}.$$

Докажем теперь, что всякое число из этого списка единственным образом представляется как сумма подчёркнутых чисел с рациональными коэффициентами. То, что неподчёркнутое число представляется нужным образом через подчёркнутые --- это по определению так. Подчёркнутое тоже представляется: оно равно само себе. Если же такое представление не единственно, то, вычитая одно представление из другого, мы видим, что сумма неких подчёркнутых чисел с рациональными коэффициентами равна 0. Но это позволяет выразить, опять-таки в виде суммы с рациональными коэффициентами, одно из подчёркнутых чисел через предыдущие подчёркнутые числа, а значит, мы зря его подчеркнули.

Все подчёркнутые числа, начиная с третьего, и будут искомым набором $e_1$,~$e_2$,~\ldots,~$e_k$.\9% Отметим, что при построении базиса мы нигде не пользовались тем, что числа, которые мы сначала записали в строку, являются длинами сторон прямоугольников какого-то разрезания, --- базис можно строить по любому набору чисел.\9% Его главное свойство сформулировано в виде задачи~\ref{l0}.% В следующих трёх задачах через \fo{(e_1,e_2,\ldots,e_k)} будем обозначать такой базис.

%\index{Базис}
Зафиксируем набор чисел $s_0$, $t_0$, $e_1$, $e_2$, \ldots, $e_n$ из леммы~\ref{l0}. Он называется \textit{базисом}.

%\problem Докажите, что подчёркнутые по такому алгоритму числа всегда образуют базис, какой бы набор чисел \fo{\lbrace s_1, s_2,\ldots,s_n\rbrace} мы не подали на вход этого алгоритма.

%\problem\label{l0} Докажите, что стороны каждого из прямоугольников \fo{s_0\times t_0}, \fo{s_1\times t_1}, \fo{s_2\times t_2}, \ldots, \fo{s_n\times t_n} (то есть каждое из чисел, по которым мы строили базис) можно единственным образом записать в виде \fo{a_1e_1+a_2e_2+\ldots+a_ke_k}.

%\2Забудем\3 старое определение \fo{x\mbox{-площади}}. Пусть у нас есть какое-то разрезание прямоугольника на прямоугольники. Построим по этому разрезанию базис \fo{(e_1,e_2,\ldots,e_k)}. Отныне 
\opred Пусть $y$ --- некоторое действительное число. Назовём ${y\mbox{-площадью}}$ прямоугольника со сторонами $$as_0+bt_0+a_1e_1+a_2e_2+\ldots+a_ne_n\,\,\,\mbox{и}\,\,\,cs_0+dt_0+c_1e_1+c_2e_2+\ldots+c_ne_n$$ число ${(a+by)(c+dy)}$. 

Обратите внимание, что при хороших несоизмеримых $s_0$ и $t_0$ и ${y=x}$ это определение не всегда эквивалентно определению ${x\mbox{-площади}}$ выше!%, если ${n>1}$, и число $a_1b_1$, если ${n=1}$.

%Ясно, что такая \fo{x\mbox{-площадь}} зависит от выбора базиса, но, ввиду громоздкости словосочетания \2\fo{x\mbox{-площадь}} прямоугольника \fo{A\times B} в базисе \fo{(e_0,e_1,\ldots,e_k)}\3, в дальнейшем будем её называть просто \fo{x\mbox{-площадью}}.
\problem Вычислите \fo{y\mbox{-площадь}} разрезаемого прямоугольника \fo{s_0\times t_0}. Является ли она неотрицательной при всех~$y$?

%\problem Докажите, что это определение корректно, то есть, что после фиксации базиса \fo{x\mbox{-площадь}} однозначно определена для любого прямоугольника, стороны которого можно записать в этом базисе.

%\refstepcounter{prob}\textbf{Задача \arabic{prob} (аддитивность \textit{x}-площади).}\label{z13} Пусть \fo{(e_0,e_1,\ldots,e_k)} --- базис, а прямоугольник \fo{1\times\left(p_0e_0+p_1e_1+\ldots+p_ne_n\right)}, где $p_0$,~$p_1$,~\ldots,~$p_n$ рациональны, разрезан\\
%а) на два;\\
%б) на произвольное число\\
%прямоугольников с длинами сторон вида \fo{c_0e_0+c_1e_1+\ldots+c_ne_n}, где $c_0$,~$c_1$,~\ldots,~$c_n$ рациональны. 
\lemma\label{l5} \textit{Любой квадрат в разрезании прямоугольника \fo{s_0\times t_0} имеет неотрицательную \fo{y\mbox{-площадь}}.}

\textbf{Доказательство.} \8Сторона любого квадрата в разрезании прямоугольника \fo{s_0\times t_0} записывается в виде \fo{as_0+bt_0+a_1e_1+a_2e_2+\ldots+a_ne_n}. Тогда его \fo{y\mbox{-площадь}} равна \fo{(a+by)^2}, что неотрицательно при любом \fo{y\in\R}.\9

\lemma\label{l6} \textit{Для любого $y$ \fo{y\mbox{-площадь}} разрезаемого прямоугольника \fo{s_0\times t_0} равна сумме \fo{y\mbox{-площадей}} прямоугольников, на которые он разрезан.}

Доказательство леммы~\ref{l6} дословно повторяет доказательство леммы~\ref{pr4} с той лишь разницей, что нужно заменить в нём \fo{x\mbox{-площадь}} на \fo{y\mbox{-площадь}}, а также заменить все числа вида \fo{a+b\sqrt{2}} на соответствующие числа вида \fo{as_0+bt_0+a_1e_1+a_2e_2+\ldots+a_ne_n} (в том числе, и на рис.~\ref{rrr6}).

%\problem Рассмотрим прямоугольник \fo{1\times\left(c_0e_0+c_1e_1+\ldots+c_ne_n\right)}, где\\1) числа $c_0$,~$c_1$,~\ldots,~$c_n$ рациональны;\\2) не существует таких рациональных $k_0$,~$k_1$,~\ldots,~$k_n$, которые не равны нулю одновременно, что \fo{k_0e_0+k_1e_1+\ldots+k_ne_n=0}.\\Докажите, что этот прямоугольник нельзя разрезать на квадраты со сторонами вида \fo{a_0e_0+a_1e_1+\ldots+a_ne_n}, где $a_0$,~$a_1$,~\ldots,~$a_n$ рациональны.

%\glava{Доказательство теоремы Дена в общем случае}

\textbf{Доказательство теоремы~\ref{D}}. \8Пусть прямоугольник \fo{s_0\times t_0} разрезан на квадраты, причём $s_0$ и $t_0$ несоизмеримы. По определению его \fo{y\mbox{-площадь}} равна $y$. Это число отрицательно при \fo{y<0}. В то же время, \fo{y\mbox{-площадь}} любого квадрата неотрицательна по лемме~\ref{l5}. Получили противоречие с леммой~\ref{l6}: сумма неотрицательных чисел не может равняться отрицательному.\9

В заключение отметим, что \fo{y\mbox{-площадь}} --- это не просто олимпиадный трюк, а частный случай \textit{\2меры\3} --- важного понятия в математике.

\refstepcounter{glav}\begin{center}\textbf{\arabic{glav}. Решения некоторых задач}\footnote{Написаны Л.~Алиевой.}\end{center}

\textbf{Решение для третьего вопроса задачи~\ref{z1}.} \textit{Ответ:} нет.

\8Предположим, что прямоугольник \fo{1\times\sqrt{2}} разрезан на $n$~квадратов со сторонами \fo{a_1+b_1\sqrt{2}}, \fo{a_2+b_2\sqrt{2}}, \ldots, \fo{a_n+b_n\sqrt{2}}. Площади квадратов соответственно равны $$a_1^2+2a_1b_1\sqrt{2}+2b_1^2,\,\,\,a_2^2+2a_2b_2\sqrt{2}+2b_2^2,\,\,\,\ldots,\,\,\,a_n^2+2a_nb_n\sqrt{2}+2b_n^2.$$ Введём обозначения: $$A:=a_1^2+a_2^2+\ldots+a_n^2,\,\,\,B:=b_1^2+b_2^2+\ldots+b_n^2,\,\,\,C:=a_1b_1+a_2b_2+\ldots+a_nb_n.$$ Тогда $A$, $B$, $C$ --- рациональные числа. Площадь прямоугольника равна сумме площадей квадратов, значит, \fo{\sqrt{2}=A+2C\sqrt{2}+B}. Получаем, что \fo{\sqrt{2}\cdot(1-2C)=A+B}. Если \fo{C=1/2}, то \fo{A=B=0}, значит, \fo{a_1=\ldots=a_n=b_1=\ldots=b_n=0}, что невозможно. Следовательно, обе части уравнения можно разделить на \fo{1-2C}. Таким образом, \fo{\sqrt{2}=(A+B)/(1-2C)}. Получаем противоречие, так как иррациональное число не может быть равно частному двух рациональных. Здесь мы вновь опирались на аддитивность площади.\9

\textbf{Решение задачи~\ref{z4}.} \textit{Ответ:} все прямоугольники с рациональным отношением сторон.

\8Рассмотрим прямоугольник \fo{\left(a+b\sqrt{2}\right)\times\left(c+d\sqrt{2}\right)}. Его \fo{x\mbox{-площадь}} равна \fo{(a+bx)(c+dx)}. Мы хотим, чтобы она была неотрицательна, то есть $$(a+bx)(c+dx)\br0\Leftrightarrow bdx^2+(bc+ad)x+ac\br0.$$ Решением данного неравенства могут быть два луча, отрезок, вся действительная ось~$Ox$ или точка. Нам надо выяснить, когда данный многочлен неотрицателен при любых $x$. То есть в каком случае решением неравенства является действительная ось~$Ox$. Если многочлен имеет два различных корня, то решением является либо отрезок, либо два луча. Следовательно, нас интересует только тот случай, когда корни многочлена совпадают, то есть \fo{x=-a/b=-c/d}. Тогда \fo{c+dx} можно представить как \fo{k(a+bx)}, где $$k=c/a=d/b=\frac{c+dx}{a+bx},\,\,\,k\in\q.$$ Заметим, что при постоянных $a$, $b$, $c$, $d$ и изменяющихся $x$ коэффициент $k$ --- константа. Однако, в случае отрицательного~$k$, решением неравенства является точка. Докажем, что данный коэффициент всегда положителен. Допустим, что \fo{k\mr0}. Возьмём \fo{x=\sqrt{2}}. Получаем, что ${\frac{c+d\sqrt{2}}{a+b\sqrt{2}}\mr0}$, то есть что отношение сторон прямоугольника неположительно, чего быть не может. Следовательно, $k$ всегда неотрицателен, а решением данного неравенства при совпадении корней всегда является вся действительная ось~$Ox$. Получаем, что \fo{x\mbox{-площадь}} прямоугольника с хорошими сторонами неотрицательна при всех $x$ тогда и только тогда, когда отношение сторон этого прямоугольника рационально.\9

\end{document}